\documentclass[11pt]{amsart}
\usepackage{latexsym,amssymb}
\usepackage{amsmath}
\usepackage{amscd}
\newcommand{\llar}{-\kern-5pt-\kern-5pt\longrightarrow}

\newtheorem{Theorem}{Theorem}[section]

\newtheorem{Corollary}[Theorem]{Corollary}

\def\Ass{\mbox{\rm Ass}}

\def\depth{\mbox{\rm depth}}
\def\ds{\displaystyle}

\def\gr{\mbox{\rm gr}}

\def\Hom{\mbox{\rm Hom}}

\def\im{\mbox{\rm Im}}

\def\l{\lambda}
\def\lar{\longrightarrow}

\def\m{\mathfrak{m}}

\def\n{\mathfrak{n}}

\def\q{\mathfrak{q}}

\def\QED{\hfill$\Box$}

\def\rar{\rightarrow}

\def\Tor{\mbox{\rm Tor}}

\begin{document}

\title{Negativity Conjecture for the First Hilbert Coefficient}

\author{L. Ghezzi}
\address{Department of Mathematics, New York City College of Technology-Cuny, 300 Jay Street, Brooklyn, NY 11201, U. S. A.}
\email{lghezzi@citytech.cuny.edu}
\author{S. Goto}
\address{Department of Mathematics, School of Science and Technology, Meiji University, 1-1-1 Higashi-mita, Tama-ku, Kawasaki 214-8571, Japan}
\email{goto@math.meiji.ac.jp}
\author{J. Hong}
\address{Department of Mathematics, Southern Connecticut State University, 501 Crescent Street, New Haven, CT 06515-1533, U. S. A.}
\email{hongj2@southernct.edu}
\author{K. Ozeki}
\address{Department of Mathematics, School of Science and Technology, Meiji University, 1-1-1 Higashi-mita, Tama-ku, Kawasaki 214-8571, Japan}
\email{kozeki@math.meiji.ac.jp}
\author{T.T. Phuong}
\address{Department of Information Technology and Applied Mathematics,
Ton Duc Thang University, 98 Ngo Tat To Street, Ward 19, Binh Thanh District,
Ho Chi Minh City, Vietnam}
\email{sugarphuong@gmail.com}
\author{W. V. Vasconcelos}
\address{Department of Mathematics, Rutgers University, 110 Frelinghuysen Rd, Piscataway, NJ 08854-8019, U. S. A.}
\email{vasconce@math.rutgers.edu}

\thanks{{AMS 2000 {\em Mathematics Subject Classification:}
13H10, 13H15, 13A30.}\\The first author is partially supported by a grant from the City University of New York PSC-CUNY Research Award Program-40. The second author is partially supported by Grant-in-Aid for Scientific Researches (C) in Japan (19540054).
The fourth author is supported by a grant from MIMS (Meiji Institute for Advanced Study of Mathematical Sciences).
The fifth author is supported by JSPS Ronpaku (Dissertation of PhD) Program.
The last author is partially supported by the NSF}


\begin{abstract}
This gives an alternate proof of the result of \cite[Theorem 2.1]{chern3}: The first Hilbert coefficient of parameter ideals in an unmixed Noetherian local ring is always negative unless the ring is Cohen--Macaulay.
\end{abstract}

\maketitle

\section{Introduction}

Let $R$ be a Noetherian local ring with the maximal ideal $\m$ of dimension $d >0$.
Let $I$ be an $\m$--primary $R$--ideal.  For sufficiently large $n$, the length $\l(R/I^{n+1})$ is of polynomial type :
\[ P_{I}(n)=\sum_{i=0}^{d} (-1)^i e_i(I) {{n+d-i}\choose{d-i}}. \]
The integers $e_i(I)$'s are called the {\em Hilbert coefficients} of $I$.
The first Hilbert coefficient $e_1(Q)$ of a parameter ideal $Q$ codes
structural information about the ring $R$ itself.
In response to a question in \cite{chern1}, the following was settled  in \cite{chern3}.

\begin{Theorem}\label{e1neg}{\rm (\cite[Theorem 2.1]{chern3})}
An unmixed Noetherian local ring $R$ is not Cohen-Macaulay if and only if $e_1(Q)<0$ for
a parameter ideal $Q$.
 \end{Theorem}

\bigskip

Meanwhile, for any parameter ideal $Q$ of $R$, it was proved that $e_1(Q)\leq 0$ (\cite[Corollary 2.5]{chern3}, \cite{MV}). Hence the above theorem can be rephrased as follows :

\begin{Corollary}\label{e1zero}
An unmixed Noetherian local ring $R$ is Cohen--Macaulay if and only if $e_1(Q)=0$ for some parameter ideal $Q$ of $R$.
\end{Corollary}

In the following section, we give an alternate proof.

\bigskip

\section{The Proof}

\noindent {\bf Proof of Theorem~\ref{e1neg}} \quad We use a setup developed in \cite{chern2}. It is enough to show that if $R$ is not Cohen--Macaulay, then $e_1(Q) <0$.
We may assume that the residue field is infinite.

\medskip

By passing to $\m$--adic completion $\widehat{R}$, we may also assume that $R$ is complete. Then there exists a Gorenstein local ring $(S, \n)$ of dimension $d=\dim(R)$ such that $R$ is a homomorphic image of $S$. This means that there exists a canonical module $\omega_R = \Hom_S(R, S)$. Consider the natural homomorphism
\[ \varphi : R \lar \Hom_S(\omega_R, \; \omega_R) \simeq \Hom_S(\omega_R,\; S). \]
Because $R$ is unmixed, this map $\varphi$ is injective (\cite[1.11.1]{Aoy}). Moreover $H^1_{\n}(R)$ has finite length.
Indeed, let $A=\Hom_S(\omega_R, \; \omega_R)$. Then by applying local cohomology to $0 \rar R \rar A \rar D \rar 0$, we obtain $H^{1}_{\n}(R) \simeq H^{0}_{\n}(D)$ since
$\depth(A) \geq 2$ (\cite{GNa}, \cite{HU}).

\bigskip

By dualizing $S^n \rar \omega_R \rar 0$ into $S$, we obtain another injective map
\[ 0 \lar \Hom_S(\omega_R,\; S) \lar S^n.  \]
Composing these two maps, we obtain an embedding
${\ds  R \hookrightarrow S^n}$.

\bigskip

Let $Q$ be a parameter ideal of $R$. Then there exists a parameter ideal $\q$ of $S$ such that $\q R =Q$ (\cite[Lemma 3.1]{chern2}). Therefore the associated graded ring of $Q$ is isomorphic to the associated graded module of $\q$ with respect to the $S$--module $R$ :
\[ \gr_Q(R) \simeq \gr_{\q}(R), \]
which implies that
\[ e_1(Q) = e_1(\q, R),  \]
where $e_1(\q, R)$ denotes the first Hilbert coefficient of $\q$ with respect to $S$--module $R$.

\bigskip

Consider the exact sequence of $S$--modules:
\[ 0 \lar R \lar S^n \lar C \lar 0.   \]
Let $y$ be a superficial element for $\q$ with respect to $R$ such that $y$ is a part of minimal generating set of $\q$ and that $y \not\in \Ass_S(C) \setminus \{ \n \}$.
By tensoring the exact sequence of $S$--modules with $S/(y)$, we get
\[ 0 \lar T =\Tor_1^S(S/yS, C) \lar  R/yR \stackrel{\zeta}{\lar} S^n/y S^n \lar C/yC \lar 0.   \]
Let $R\,'= R/yR$ and $S\,'=\im(\zeta)$ and consider the short exact sequence :
\[ 0 \lar T \lar R\,'  \lar S\,' \lar 0. \]
Then either $T=0$ or $T$ has finite length $\l(T) < \infty$.

\bigskip

\noindent Now we use induction on $d=\dim(R)$ to show that if $R$ is not Cohen--Macaulay, then $e_1(\q, S) <0$.

\bigskip

Let $d=2$ and $\q=(y, z)$. Then $T \neq 0$ so that $\l(T) < \infty$. Applying the Snake Lemma to
\[
\begin{CD}
 0 @>>> T \cap z^n R\,' @>>> z^n R\,' @>>> z^n S\,' @>>> 0 \\
  & & @VVV  @VVV  @VVV \\
0 @>>> T @>>> R\,' @>>> S\,' @>>> 0
\end{CD}
\] we get, for sufficiently large $n$,
\[  \l(R\,'/z^n R\,') = \l(T) + \l(S\,'/z^n S\,').  \]
Computing the Hilbert polynomials, we have
\[ e_1(\q/y, R/yR) = -\l(T) <0 \]
so that
\[ e_1(\q, R) = e_1(\q/y, R/yR) - \l(0:_R y) = -\l(T) - \l(0:_R y) < 0. \]

\bigskip

Now suppose that $d \geq 3$. From the exact sequence
\[ 0 \lar T \lar R\,'=R/yR  \lar S\,' \lar 0, \]
we have
\[ e_1(\q, R) = e_1(\q/(y) , R/yR ) = e_1(\q/(y) , S\,').  \]
By an induction argument, it is enough to show that $S\,'$ is not Cohen--Macaulay since $\dim(S/yS)=d-1$.

\bigskip

Suppose that $S\,'$ is Cohen--Macaulay. Let $\n$ be the maximal ideal of $S/yS$. From the exact sequence
\[ 0 \lar T \lar R\,'=R/yR  \lar S\,' \lar 0, \]
we obtain the long exact sequence:
\[
0 \rar H_{\n}^0(T) \rar H_{\n}^0(R\,') \rar H_{\n}^0(S\,') \rar H_{\n}^1(T) \rar H_{\n}^1(R\,') \rar H_{\n}^1(S\,').
\]By the assumption that  $S\,'$ is Cohen--Macaulay of dimension $d-1 \geq 2$ and the fact that $T$ is a torsion module, we get

\[
0 \rar T \stackrel{\simeq}{\lar} H_{\n}^0(R\,') \rar 0 \rar 0 \rar H_{\n}^1(R\,') \rar 0.
\]


\noindent We may assume that $y$ is a nonzerodivisor on $R$. From the exact sequence
\[0 \lar R \stackrel{\cdot y}{\lar} R \lar R/yR \lar 0, \]
we obtain the following exact sequence:
\[
0  \lar {T \simeq H_{\n}^0(R\,')} \lar H_{\n}^1(R) \stackrel{\cdot y}{\lar} H_{\n}^1(R) \lar {H_{\n}^1(R\,')=0}.
\] Since $H_{\n}^1(R) $ is finitely generated and ${\ds H_{\n}^1(R)=yH_{\n}^1(R)} $, we have ${\ds H_{\n}^1(R)=0}$. This means that
$T=0$. Therefore
\[ {\ds 0 \rar T=0 \rar R/yR \stackrel{\simeq}{\lar} S\,' \rar 0.} \]
Since $S\,'$ is Cohen--Macaulay, $R\,'=R/yR$ is Cohen--Macaulay.
Since $y$ is regular on $R$, $R$ is Cohen--Macaulay, which is a contradiction. \QED

\end{document}